# Knuth-Robinson-Schensted correspondence and Weak Polynomial Identities of $M_{1,1}(E)$ *

Onofrio Mario Di Vincenzo, Roberto La Scala

**Abstract**

In this paper it is proved that the ideal $I_w$ of the weak polynomial identities of the superalgebra $M_{1,1}(E)$ is generated by the proper polynomials $[x_1, x_2, x_3]$ and $[x_2, x_1][x_3, x_1][x_4, x_1]$. This is proved for any infinite field $F$ of characteristic different from 2. Precisely, if $B$ is the subalgebra of the proper polynomials of $F\langle X\rangle$, we determine a basis and the dimension of any multihomogeneous component of the quotient algebra $B/B \cap I_w$. We compute also the Hilbert series of this algebra. One of the main tools of the paper is a variant we found of the Knuth-Robinson-Schensted correspondence defined for single semistandard tableaux of double shape.

## 1 Introduction

In the structure theory of the varieties of associative algebras as developed by Kemer [13] a fundamental role is played, together with the usual matrix algebras over a field, by the superalgebras $M_{k,l}(E)$ of matrices with entries in the Grassmann algebra $E$. A main purpose in the study of $M_{k,l}(E)$ is to find out bases for ideals of polynomial identities satisfied by such algebras. In particular, Razmyslov [18, 19, 20] introduced the notion of "weak polynomial identity" for both the algebras $M_n(F)$ and $M_{k,l}(E)$, and explained how these identities are correlated with central polynomials and identities in the traces.

For fields of characteristic zero, Razmyslov [18] has found finite bases for the polynomial identities of $M_2(F)$ and $sl_2(F)$ (the Lie algebra of traceless matrices). Consequently, Drensky and Filippov [7, 10] described minimal bases for the identities of such algebras. For infinite fields of characteristic different from 2, Koshlukov in [15] has computed a basis for the weak polynomial identities of the pair $(M_2(F), sl_2(F))$ and successively in [16] described a finite basis for the ordinary polynomial identities of $M_2(F)$.

*Partially supported by Università di Bari and COFIN-MIUR 2001 "Algebre con identità polinomiali e metodi combinatori"

In the characteristic zero case, it is well known (see [13]) that $M_{1,1}(E)$ and $E \otimes E$ satisfy the same ordinary polynomial identities. Popov in [17] has computed a basis for such identities and Kemer in [12] has considered the weak polynomial identities of $E \otimes E$ and related them to the ordinary ones. The ideal of $\mathbb{Z}_2$-graded identities of the superalgebra $M_{1,1}(E)$ has been studied in [6]. A basis for the weak polynomial identities of $M_{1,1}(E)$ has been found, still in characteristic zero, in [5].

For infinite fields of characteristic different from 2, Azevedo and Koshlukov in [1] have found bases for the $\mathbb{Z}_2$-graded identities of the superalgebras $M_{1,1}(E)$ and $E \otimes E$. Finally, under the same assumption for the base field, in the present paper we study the weak polynomial identities of $M_{1,1}(E)$.

More precisely, we prove that in the free algebra $F\langle X \rangle$ the ideal $I_w$ of such identities is generated by the polynomials $[x_1, x_2, x_3]$ and $[x_2, x_1][x_3, x_1][x_4, x_1]$. If $B$ is the subalgebra of the proper polynomials of $F\langle X \rangle$, we describe the multihomogeneous components of the quotient algebra $B(I_w) = B/B \cap I_w$. We compute bases and dimensions of such components, and also the Hilbert series of the algebra $B(I_w)$. Our approach is essentially combinatorial. One of the techniques is the study of the invariant ring of the orthogonal group as developed by De Concini and Procesi [4] in a characteristic free way. Another tool is a variant we found of the Knuth-Robinson-Schensted correspondence which is defined for single semistandard tableaux of double shape.

## 2 Basics

Let $F$ be any field and $F\langle X \rangle$ the free associative algebra generated by a countable set of variables $X = \{x_1, x_2, \dots\}$. If $R$ is an associative algebra and $W \subset R$ is a vector space, then the polynomial $f(x_1, \dots, x_n) \in F\langle X \rangle$ is called *weak polynomial identity* for the pair $(R, W)$ if $f(w_1, \dots, w_n) = 0$, for all the elements $w_1, \dots, w_n \in W$. The set of all weak polynomial identities is an ideal $I_w = I(R, W)$ of the ring $F\langle X \rangle$. It is well known that for an suitable description of $I_w$, it is convenient to determine endomorphisms of $F\langle X \rangle$ which stabilize $I_w$, that is to establish rules that allow to take consequences from any set of weak polynomial identities. More precisely, let $\Omega$ be a non-empty subset of $F\langle X \rangle$ such that $\omega(w_1, \dots, w_n) \in W$, for all choices of $\omega \in \Omega$ and $w_1, \dots, w_n \in W$. Then, the ideal $I_w$ is stable under the endomorphisms of the algebra $F\langle X \rangle$ corresponding to polynomials of $\Omega$. In general, any ideal $I \subset F\langle X \rangle$ which verifies such property is called $\Omega$-*stable* or simply $\Omega$-*ideal*. Let now $G$ be a non-empty subset of $F\langle X \rangle$. A polynomial $f(x_1, \dots, x_n)$ is said to be an $\Omega$-*consequence* of $G$ is $f$ belongs to the minimal $\Omega$-stable ideal $I \supset G$. In this case, $G$ is called an $\Omega$-*generating set* of the ideal $I$. Depending on the properties of the space $W$, it is possible to choose the set $\Omega$ in different ways. The simplest choice is to put $W = R$ and $\Omega = F\langle X \rangle$, so that $I_w = T(R)$ is the $T$-ideal of ordinary polynomial identities of $R$.

If we assume that the field $F$ is infinite, it is well-known that any (weak) polynomial identity $f = 0$ is equivalent to the collection of identities given by the multihomogeneous components of the polynomial $f$. We may therefore reduce the study of (weak) polynomial identities to the multihomogeneous case. Let $B$ denote the subalgebra of $F\langle X\rangle$ generated by all the commutators $[x_{i_1},\ldots,x_{i_l}]$ of length $l \geq 2$. For a finite number of variables, we put $B_k = B \cap F\langle x_1,\ldots,x_k\rangle$. The elements of $B$ are called *proper polynomials*. It is well-known that any T-ideal $I$ of $F\langle X\rangle$ is generated by $B \cap I$. Assume now that $\Omega$ contains the subspace of $F\langle X\rangle$ spanned by the set $X \cup \{1\}$ (hence $1 \in W$). By the same argument used for the T-ideals, it is possibile to show that any $\Omega$-ideal $I$ is $\Omega$-generated by the set $B \cap I$ (see [5], [8]).

Let $I$ be any $\Omega$-ideal. We define the quotient algebra $B(I) = B/B \cap I$ and similarly $B_k(I) = B_k/B_k \cap I$, for any $k > 0$. If $(n_i) = (n_1,\ldots,n_k)$ is a vector of integers $n_i \geq 0$, for the multigraded algebra $B$ we denote by $B^{(n_i)}$ the multihomogeneous component of $B$ of *multidegree* $(n_i)$. Under the assumption that the subspace spanned by $X \cup \{1\}$ is a subset of $\Omega$, we have that the $\Omega$-ideal $I$ is multihomogeneous and the quotient algebra $B(I)$ multigraded. Then, the components of $B(I)$ are:

$$B^{(n_i)}(I) = \frac{B^{(n_i)}}{B^{(n_i)} \cap I}$$

Hence, for studying the $\Omega$-ideal $I$ it is sufficient to describe the vector spaces $B^{(n_i)}(I)$, for all the multidegrees $(n_i)$.

Let $\bar{F}$ denote the algebraic closure of the field $F$ and consider the algebra $R_{\bar{F}} = R \otimes_F \bar{F}$ and the subspace $W_{\bar{F}} = W \otimes_F \bar{F}$ obtained by extending the scalars from $F$ to $\bar{F}$. Since $F$ is infinite, we have:

$$I(R_{\bar{F}}, W_{\bar{F}}) = I(R, W) \otimes_F \bar{F}$$

that is the algebras $R, R_{\bar{F}}$ satisfy the same (proper, weak) polynomial identities over $F$.

Finally, let $F$ be any infinite field of characteristic different from 2. We denote by $E = E_0 \oplus E_1$ the Grassmann algebra generated by a vector space of countable dimension over $F$, and by $E_0, E_1$ the $\mathbb{Z}_2$-homogeneous components of $E$. We define $R$ the following matrix superalgebra:

$$R = M_{1,1}(E) = \{\begin{pmatrix} a & b \\ c & d \end{pmatrix} \mid a,d \in E_0,\ b,c \in E_1\}$$

For any matrix $A \in R$, a *supertrace* is defined in the following way:

$$\text{str}(A) = a - d$$

This supertrace verifies the usual properties of traces. Denote by $W$ the vector space given by all the matrices of $R$ with supertrace equal to zero. We define then $I_w = I(R, W)$ the ideal of $F\langle X\rangle$ of the weak polynomial identities for the pair $(R, W)$. Let now $\Omega$ be the subspace of $F\langle X\rangle$ spanned by $X \cup \{1\}$, that is its elements are all the linear polynomials. The main purpose of the present paper is to compute an $\Omega$-generating set for the ideal $I_w \subset F\langle X\rangle$. From what we have observed before, we may assume, without loss of generality, that the field $F$ is algebraically closed.

# 3  Identities and spanning

We note immediately that $c_3 := [x_1, x_2, x_3] \in B$ is a weak polynomial identity for the algebra $R = M_{1,1}(E)$. Precisely, for any $w_1, w_2 \in W$ the commutator $[w_1, w_2]$ is in the center of $R$. Since $c_3 \in I_w$ we have that $B^{(n_i)}(I_w) = 0$ if the *total degree* $\sum n_i$ is odd. Hence, in what follows we shall assume that $\sum n_i = 2m$. Let now $S$ be any two-rowed array of type:

$$S = \left[ \begin{array}{cccc} a_1 & a_2 & \ldots & a_m \\ b_1 & b_2 & \ldots & b_m \end{array} \right]$$

with $a_i, b_i > 0$ integers. Denote by $f_S$ the following proper polynomial:

$$f_S = [x_{a_1}, x_{b_1}][x_{a_2}, x_{b_2}] \cdots [x_{a_m}, x_{b_m}]$$

Note that if $f_S$ has multidegree $(n_i) = (n_1, \ldots, n_k)$ then the array $S$ has *content* $(n_i)$, i.e. in the multiset $\{a_1, b_1, \ldots, a_m, b_m\}$ the integer 1 occurs $n_1$ times, 2 occurs $n_2$ times, etc, and one has $\sum n_i = 2m$.

Since $[x_i, x_j] = -[x_j, x_i]$ and $c_3 \in I_w$, we have clearly that any polynomial $f_S \neq 0$ can be rewritten modulo $I_w$ so that the indices $a_i, b_i$ satisfy the following conditions:

s$_1$) $a_i > b_i$, for all $i = 1, 2, \ldots, m$

s$_2$) $(a_i, b_i) \leq (a_{i+1}, b_{i+1})$ w.r.t. (left) lexicographic ordering

For this reason, we call *c-array* or *array of commutators* any two-rowed array which verifies (s$_1$) and (s$_2$). Formally, an array $S$ can be transformed into a c-array by applying consecutively the following rules:

R$_1$) in any column $(a_i, b_i)$ of $S$, swap $a_i$ and $b_i$ if $a_i < b_i$

R$_2$) sort the columns of $S$ with respect to the lexicographic ordering

It is easy to prove that the algebra $R$ satisfies another weak polynomial identity.

**Proposition 3.1** *The proper polynomial:*

$$p := [x_2, x_1][x_3, x_1][x_4, x_1]$$

*is an element of the ideal $I_w$.*

Denote now by $I \subset F\langle X \rangle$ the ideal $\Omega$-generated by the polynomials $c_3$ and $p$. We shall prove that $I_w = I$. From what we have observed for the identity $c_3$, we get immediately:

**Lemma 3.2** *A generating set for the vector space $B^{(n_i)}(I)$ is given by the proper polynomials $f_S$, where $S$ ranges over the c-arrays of content $(n_i)$.*

Actually, the rewriting action of the identity $p$ gives rise to a smaller generating set. Note that by the identities $[x_i, x_j] = -[x_j, x_i]$, $c_3$ and $p$, we have $B^{(n_i)}(I_w) = B^{(n_i)}(I) = 0$ if $n_i > 2$, for some $i$. Then, we say that a c-array $S$ is *normal* if it verifies also the following conditions:

s$_3$) $n_i \leq 2$, for all $i = 1, 2, \ldots, k$

s$_4$) there are no indices $r < s < t$ such that $b_r \leq b_s \leq b_t$

Denote by $N^{(n_i)}$ the set of normal c-arrays of content $(n_i)$. Note that the set of all two-rowed arrays is totally ordered in the following way:

$$\begin{bmatrix} a_1 & \ldots & a_m \\ b_1 & \ldots & b_m \end{bmatrix} < \begin{bmatrix} c_1 & \ldots & c_m \\ d_1 & \ldots & d_m \end{bmatrix} \qquad (1)$$

if and only if $(a_m, \ldots, a_1, b_1, \ldots, b_m) < (c_m, \ldots, c_1, d_1, \ldots, d_m)$ in the lexicographic order. Let now $S', S''$ be any pair of c-arrays. By juxtaposing these arrays and then applying the rule $R_2$ we get a unique c-array $S = S' * S''$. Clearly, the equality $f_S = f_{S'} \cdot f_{S''}$ holds modulo $I$. Moreover, it is easy to prove that for the c-arrays the ordering (1) is *compatible* with the product "$*$" that is:

$$S' < S'' \iff S * S' < S * S''$$

for any $S, S', S''$ c-arrays. We are ready to improve Lemma 3.2.

**Theorem 3.3** *A generating set for $B^{(n_i)}(I)$ is given by the proper polynomials $f_S$, for all $S \in N^{(n_i)}$.*

*Proof:* For simplicity, we call *almost-normal* the c-arrays $S$ which satisfy condition (s$_3$). By Lemma 3.2 we have that the polynomials $f_S$, where $S$ is any almost-normal c-array, define a generating set for the space $B^{(n_i)}(I)$. Then, it is sufficient to prove that if $S$ does not verify also the condition (s$_4$) then $f_S$ is a linear combination of proper polynomials corresponding to almost-normal c-arrays which are greater in the ordering (1).

Of course, any almost-normal c-array is actually normal for $m < 3$. Then, we argue for $m = 3$. Consider:

$$S = \begin{bmatrix} a_1 & a_2 & a_3 \\ b_1 & b_2 & b_3 \end{bmatrix}$$

an almost-normal c-array which does not verify condition (s$_4$) i.e. $b_1 \leq b_2 \leq b_3$. Since the $\Omega$-ideal $I$ is stable under linear substitutions of variables and the field $F$ is infinite, we have that any (partial) linearization of $p$ is still an element of $I$. If the integers $b_1, b_2, b_3$ are all distict then, modulo the ideal $I$, we have:

$$\sum_{\sigma \in \mathbb{S}_3} f_{\sigma(S)} = 0, \text{ where } \sigma(S) = \begin{bmatrix} a_1 & a_2 & a_3 \\ b_{\sigma(1)} & b_{\sigma(2)} & b_{\sigma(3)} \end{bmatrix} \qquad (2)$$

Otherwise, if $\{b_1, b_2, b_3\} = \{b, \beta\}$ (necessarily $b \neq \beta$) then we have:

$$f_{S_1} + f_{S_2} + f_{S_3} = 0 \tag{3}$$

where:

$$S_1 = \begin{bmatrix} a_1 & a_2 & a_3 \\ b & b & \beta \end{bmatrix}, S_2 = \begin{bmatrix} a_1 & a_2 & a_3 \\ b & \beta & b \end{bmatrix}, S_3 = \begin{bmatrix} a_1 & a_2 & a_3 \\ \beta & b & b \end{bmatrix}$$

and $S = S_1$ holds if $b < \beta$ or $S = S_3$ if $\beta < b$.

Since $b_1 \leq b_2 \leq b_3$, any array $S^* \neq S$ which occurs in the equations (2),(3) is strictly greater than $S$ in the ordering. If $S^*$ is not a c-array and $f_{S^*} \neq 0$, then we may define a unique c-array $\bar{S}^*$ by applying the rules $R_1, R_2$ to $S^*$. Then, we substitute in the above equations the polynomial $f_{S^*}$ with $\pm f_{\bar{S}^*}$ which is equivalent to it modulo $I$. Note that if $R_1(S^*) \neq S^*$ then in $R_1(S^*)$ we have substituted some integers of the first row of $S^*$ with integers from the second row which are strictly greater. Hence, after we have eventually sorted the columns of $R_1(S^*)$, it holds:

$$\bar{S}^* = R_2(R_1(S^*)) > S^* > S$$

If instead $R_1(S^*) = S^*$ and $R_2(S^*) \neq S^*$, then the array $S^*$ is of type:

$$S^* = \begin{bmatrix} a_1 & a_2 & a_3 \\ b_1^* & b_2^* & b_3^* \end{bmatrix}$$

where $a_i > b_i^*$ for all $i$, two of the integers $a_1 \leq a_2 \leq a_3$ are equal, say $a_i, a_{i+1}$, and $b_i^* > b_{i+1}^*$ holds. In this case, since $\mathrm{char}(F) \neq 2$ we may rewrite the equation (2) as:

$$f_S + f_{S'} + f_{S''} = 0 \tag{4}$$

where, if $a_1 = a_2$ one has:

$$S = \begin{bmatrix} a & a & \alpha \\ b_1 & b_2 & b_3 \end{bmatrix}, S' = \begin{bmatrix} a & a & \alpha \\ b_1 & b_3 & b_2 \end{bmatrix}, S'' = \begin{bmatrix} a & a & \alpha \\ b_2 & b_3 & b_1 \end{bmatrix}$$

Otherwise, if $a_2 = a_3$ it holds:

$$S = \begin{bmatrix} \alpha & a & a \\ b_1 & b_2 & b_3 \end{bmatrix}, S' = \begin{bmatrix} \alpha & a & a \\ b_2 & b_1 & b_3 \end{bmatrix}, S'' = \begin{bmatrix} \alpha & a & a \\ b_3 & b_1 & b_2 \end{bmatrix}$$

In the same way, if $a_1 = a_2$ then the equation (3) can be rewritten as:

$$f_S + 2f_{S'} = 0, \text{ with } S = \begin{bmatrix} a & a & \alpha \\ b & b & \beta \end{bmatrix}, S' = \begin{bmatrix} a & a & \alpha \\ b & \beta & b \end{bmatrix} \text{ if } b < \beta \tag{5}$$

$$2f_S + f_{S'} = 0, \text{ with } S = \begin{bmatrix} a & a & \alpha \\ \beta & b & b \end{bmatrix}, S' = \begin{bmatrix} a & a & \alpha \\ b & b & \beta \end{bmatrix} \text{ if } \beta < b \tag{6}$$

Otherwise, if $a_2 = a_3$ it holds:

$$2f_S + f_{S'} = 0, \text{ with } S = \begin{bmatrix} \alpha & a & a \\ b & b & \beta \end{bmatrix}, S' = \begin{bmatrix} \alpha & a & a \\ \beta & b & b \end{bmatrix} \text{ if } b < \beta \qquad (7)$$

$$f_S + 2f_{S'} = 0, \text{ with } S = \begin{bmatrix} \alpha & a & a \\ \beta & b & b \end{bmatrix}, S' = \begin{bmatrix} \alpha & a & a \\ b & \beta & b \end{bmatrix} \text{ if } \beta < b \qquad (8)$$

Note again that each array $S^* \neq S$ which occurs in the equations (4)–(8), verifies $S^* > S$. Moreover, if $f_{S^*} \neq 0$ then the unique c-array $\bar{S}^*$ corresponding to $S^*$ is greater than $S^*$ In fact, in these equations one has $R_2(S^*) = S^*$ whenever $R_1(S^*) = S^*$.

We argue now for any $m > 3$. If $S$ is any almost-normal c-array that does not verify condition $(s_4)$, then there exists an array $U$ formed by 3 columns of $S$ (not necessarily consecutive) which is almost-normal and does not satisfies $(s_4)$. Since $c_3 \in I$, the proper polynomial $f_S$ can be rewritten modulo $I$ as:

$$f_S = f_U \cdot f_V$$

where $V$ is an array obtained by complementing $U$ in $S$ so that $S = U * V$. We have proved for $m = 3$ that $f_U$ can be rewritten modulo $I$ as a linear combination of proper polynomials associated to almost-normal c-arrays which are greater than $U$ in the total ordering (1). Then, it happens also for $f_S$ owing to compatibility of such ordering with the product "$*$".

■

## 4 From c-arrays to d-tableaux

Let $\lambda = (\lambda_1, \ldots, \lambda_r)$ be any partition of an integer $n > 0$. We call *tableau of shape* $\lambda$ and *degree* $n$ simply any array of positive integers $T = (t_{ij})$, with indices in the ranges $1 \leq i \leq r$ and $1 \leq j \leq \lambda_i$. In particular, the c-arrays are special tableaux of shape $(m, m)$. The vector $(n_1, \ldots, n_k)$ of integers $n_i \geq 0$ is called the *content* of the tableau $T$ if 1 occurs $n_1$ times in $T$, 2 occurs $n_2$ times, etc, and we have $\sum n_i = n$. A tableau $T$ is said to be *multilinear* whenever $n_i = 1$, for any $i$.

Following the "english notation", we call a tableau $T = (t_{ij})$ *semistandard* if for any $i, j$ it holds:

1. $t_{ij} \leq t_{ij+1}$ (weakly increasing across each row)

2. $t_{ij} < t_{i+1j}$ (strictly increasing down each column)

A multilinear semistandard tableau is called simply *standard*. We define *d-tableau* or *tableau of double shape* any semistandard tableau whose shape is a partition of an even integer $2m$ of type $\lambda = (\lambda_1, \lambda_1, \ldots, \lambda_r, \lambda_r)$. For simplifying the notation, we write $\lambda = (\lambda_1^2, \ldots, \lambda_r^2)$.

We will prove that a bijection is given between the set of c-arrays of some content and the set of d-tableaux of the same content. This bijection is based on the Knuth-Robinson-Schensted correspondence (see for instance [14]) and appears as a variant to the "english notation" of the correspondence found by Conca [3] for double tableaux defined semistandard in the "french notation" (strictly increasing by rows and weakly by columns). Note that the two correspondence are the same just in the multilinear case.

Since the bijection is based on KRS-correspondence, we have to introduce the fundamental algorithms of insertion and deletion by rows that we denote as Insert and Delete. The inputs of Insert are a semistandard tableau $T = (t_{ij})$ of shape $(\lambda_1, \ldots, \lambda_r)$ and an integer $x$. The outputs are a row index $1 \leq i \leq r+1$ and a semistandard tableau $T'$ of shape $(\lambda_1, \ldots, \lambda_{i-1}, \lambda_i + 1, \lambda_{i+1}, \ldots, \lambda_r)$ if $i \leq r$, or $(\lambda_1, \ldots, \lambda_r, 1)$ otherwise. In the procedure, we denote by $R_i$ the $i$-th row of $T$.

**Procedure 4.1** Insert$(T, x)$

```
add an empty row R_{r+1} to T
i := 1
J := { j | t_{ij} ∈ R_i, t_{ij} > x }
while J ≠ ∅ do
   k := min(J)
   replace the entry t_{ik} by x in R_i
   x := t_{ik}
   i := i + 1
   J := { j | t_{ij} ∈ R_i, t_{ij} > x }
append x at the end of the row R_i
if i ≤ r then delete the empty row R_{r+1} from T
RETURN(T, i)
```

The inputs of Delete are given by a semistandard tableau $T$ of shape $(\lambda_1, \ldots, \lambda_r)$ and a row index $1 \leq i \leq r$ such that $\lambda_i > \lambda_{i+1}$. The outputs are a semistandard tableau $T'$ of shape $(\lambda_1, \ldots, \lambda_{i-1}, \lambda_i - 1, \lambda_{i+1}, \ldots, \lambda_r)$ and an entry $x$ of the tableau $T$.

**Procedure 4.2** Delete$(T, i)$

```
add an empty row R_0 to T
x := the last entry of the row R_i
delete x from R_i
h := i - 1
J := { j | t_{hj} ∈ R_h, t_{hj} < x }
while J ≠ ∅ do
   k := max(J)
   replace the entry t_{hk} by x in R_h
   x := t_{hk}
   h := h - 1
   J := { j | t_{hj} ∈ R_h, t_{hj} < x }
delete the empty row R_0 from T
```

```
if R_r is also an empty row then delete it
RETURN(T, x)
```

The correctness of these algorithms is clear and it holds:

$$\texttt{Insert}(T, x) = (T', i) \iff \texttt{Delete}(T', i) = (T, x)$$

The following basic result is due to Knuth [14].

**Theorem 4.3 (Row Bumping Lemma)** *Let $T$ be a semistandard tableau and $x, y$ integers. Put $(T', i) = \texttt{Insert}(T, x)$ and $(T'', j) = \texttt{Insert}(T', y)$. Moreover, denote by $h, k$ respectively the lengths of the rows $i, j$ of the semistandard tableaux $T', T''$. It holds:*

  *i) if $x \leq y$ then $i \geq j$ and $h < k$*

  *ii) if $x > y$ then $i > j$ and $h \geq k$*

For the proofs and further details we refer the reader to the Fulton's book [11].

Following Conca [3], we define a one-to-one correspondence between the c-arrays and the d-tableaux of some fixed content by means of the procedure `Insert` and `Delete`. Such correspondence is accomplished by two algorithms that we call `Carray2Dtableau` and `Dtableau2Carray`. The procedure `Carray2Dtableau` gets as input a c-array of type:

$$S = \left[ \begin{array}{cccc} a_1 & a_2 & \ldots & a_m \\ b_1 & b_2 & \ldots & b_m \end{array} \right]$$

Recall that we have $a_k > b_k$ and $(a_k, b_k) \leq (a_{k+1}, b_{k+1})$ in the lexicographic ordering. The output is given by a d-tableau $T$ of shape $\lambda \vdash 2m$ and entries in the multiset $\{a_1, b_1, \ldots, a_m, b_m\}$. The procedure `Dtableau2Carray` has exactly the reverse input-output.

**Procedure 4.4** `Carray2Dtableau(S)`

```
T := the tableau with 2m empty rows
r := 0
for k from 1 to m do
   (T, i) := Insert(T, b_k)
   append a_k at the end of the row R_{i+1}
   if i + 1 > r then r := i + 1
if r < 2m then delete the last 2m − r empty rows of T
RETURN(T)
```

**Procedure 4.5** `Dtableau2Carray(T)`

```
for k from m to 1 do
   x := max(T)
   t_ij := the occurrence of x with maximal column j
   (T,y) := Delete(T,i)
   a_k := x; b_k := y
   delete the entry t_{i-1 j} from R_{i-1}
   if R_{i-1} is an empty row then delete it from T
S := (a_1, b_1, ..., a_m, b_m)
RETURN(S)
```

**Proposition 4.6** *The algorithms* `Carray2Dtableau` *and* `Dtableau2Carray` *are correct. Moreover, for any c-array $S$ and d-tableau $T$, we have:*

$$\texttt{Carray2Dtableau}(S) = T \Leftrightarrow \texttt{Dtableau2Carray}(T) = S$$

*Proof:* Let $T$ be a semistandard tableau of double shape $(\lambda_1^2, \ldots, \lambda_r^2)$. Put $a = \max(T)$ and denote by $j$ the maximal column where $a$ occurs. Moreover, let $a' > b'$ be integers. Put $(\hat{T}, i') = \texttt{Insert}(T, b')$ and define $j'$ the length of the row $i'$. Of course, the integer $i' = 2k - 1$ is necessarily odd. Then, define $T'$ the tableau obtained from $\hat{T}$ by appending the element $a'$ at the end of the row $i'+1$. Clearly, the tableau $T'$ has double shape $(\lambda_1^2, \ldots, \lambda_{k-1}^2, (\lambda_k+1)^2, \lambda_{k+1}^2, \ldots, \lambda_r^2)$. If we suppose that:

$$a < a' \text{ or } a = a', j < j' \tag{9}$$

we have immediately that $T'$ is semistandard. The procedure $T \overset{\texttt{Insert}}{\longrightarrow} \hat{T} \longrightarrow T'$ can be reversed in the obvious way as: $T \overset{\texttt{Delete}}{\longleftarrow} \hat{T} \longleftarrow T'$. It is also possible to obtain the tableau $T$ starting from $T'$ in an alternative way. By putting $(\check{T}, x) = \texttt{Delete}(T', i' + 1)$, under the condition (9) it is easy to verify that $x = b'$ and $a'$ occurs in the position $(i', j')$ in the tableau $\check{T}$. Then, $T$ can be obtained by $\check{T}$ simply deleting this entry. We have hence the alternative procedure: $T \longleftarrow \check{T} \overset{\texttt{Delete}}{\longleftarrow} T'$ which can be reversed in the obvious way as: $T \longrightarrow \check{T} \overset{\texttt{Insert}}{\longrightarrow} T'$. It is clear now that the algorithms `Carray2Dtableau` and `Dtableau2Carray` are one the inverse of the other if they are both correct.

Let $T \longrightarrow T' \longrightarrow T''$ be tableaux obtained by two consecutive steps of the algorithm `Carray2Dtableau` and let $(a', b'), (a'', b'')$ be the pairs of integers corresponding to $T', T''$. Moreover, denote by $j', j''$ the maximal columns in which the integers $a', a''$ occur respectively in $T', T''$. By induction, we may assume that the condition (9) is verified and therefore $T'$ is semistandard. Now, if $a' < a''$ then also $T''$ is semistandard. Otherwise, from the condition (s$_2$) that defines the c-arrays, it follows that $a' = a''$ and $b' \leq b''$. Since we may transform $T$ into $T''$ by means of:

$$T \longrightarrow \check{T} \overset{\texttt{Insert}}{\Longrightarrow} T' \overset{\texttt{Insert}}{\Longrightarrow} \hat{T}' \longrightarrow T''$$

from the Row Bumping Lemma we have that $j' < j''$ and hence $T''$ is semistandard.

We argue now for `Dtableau2Carray` i.e. for $T \longleftarrow T' \longleftarrow T''$. By induction, we have that the condition (9) is satisfied by $T' \longleftarrow T''$ that is $a' < a''$ or $a' = a'', j' < j''$. Since we may transform $T''$ into $T$ as:

$$T \longleftarrow \check{T} \overset{\text{Delete}}{\longleftarrow} T' \overset{\text{Delete}}{\longleftarrow} \hat{T}' \longleftarrow T''$$

from the Row Bumping Lemma it follows that $a' < a''$ or $a' = a'', b' \leq b''$ i.e. $(a', b') \leq (a'', b'')$ in the lexicographic ordering. ∎

The next proposition will be fundamental in proving that the normal c-arrays parametrize a linear basis of $B^{(n_i)}(I_w)$.

**Proposition 4.7** *Let $S$ be a c-array of type:*

$$S = \left[\begin{array}{cccc} a_1 & a_2 & \ldots & a_m \\ b_1 & b_2 & \ldots & b_m \end{array}\right]$$

*Put $T = \mathtt{Carray2Dtableau}(S)$ and let $l$ denote the length of the first row of the d-tableau $T$. Then $l$ is the length of the longest (weakly) increasing sequence $b_{i_1} \leq \ldots \leq b_{i_l}$ $(1 \leq i_1 < \ldots < i_l \leq m)$ one can extract from $b_1, \ldots, b_m$.*

*Proof:* The "classic" algorithm of Knuth-Robinson-Schensted transforms $S$ into a couple of semistandard tableaux $T', T''$ of the same shape and content respectively equal to the multisets $\{b_1, \ldots, b_m\}$ and $\{a_1, \ldots, a_m\}$. Then, let $R, R'$ denote the first rows respectively of $T, T'$. By the definition of `Carray2Dtableau` we have immediately that $R = R'$ and therefore the claim follows by an equivalent result for the KRS algorithm ([14], Corollary 4.2). ∎

## 5 Linear independence: the multilinear case

Denote by $P = F[U_i, V_i]$ the polynomial ring in $4m$ commuting variables $U_i, V_i$, with $1 \leq i \leq 2m$. We endow $P$ with the graduation that associates to the monomial:

$$U_1^{c_1} \ldots U_{2m}^{c_{2m}} V_1^{d_1} \ldots V_{2m}^{d_{2m}}$$

the multidegree $(c_1+d_1, \ldots, c_{2m}+d_{2m})$. Consider the graded subalgebra $\Pi \subset P$ generated by the multihomogeneous polynomials:

$$q_{ij} = U_i U_j + V_i V_j$$

for all $1 \leq i, j \leq 2m$. If $(n_i) = (n_1, \ldots, n_{2m})$ is a vector of integers $n_i \geq 0$, denote as usual by $\Pi^{(n_i)}$ the multihomogeneous component of the algebra $\Pi$ of multidegree $(n_i)$.

Let us call a tableau $T = (t_{ij})$ of double shape $(\lambda_1^2, \ldots, \lambda_r^2)$ a $d^*$-*tableau* if $T$ is semistandard according to the "french notation", that is:

1. $t_{ij} < t_{ij+1}$ (strictly increasing across each row)

2. $t_{ij} \leq t_{i+1j}$ (weakly increasing down each column)

Clearly, the d*-tableaux are equal to our d-tableaux just in the multilinear case. The next proposition follows immediately by a result of De Concini and Procesi ([4], Theorem 5.1).

**Proposition 5.1** *For the vector space $\Pi^{(n_i)}$ is given the decomposition:*

$$\Pi^{(n_i)} = \bigoplus_\lambda \Pi_\lambda^{(n_i)}$$

*where $\lambda$ ranges in the set of all partitions of type $\lambda = (2^{2p}, 1^{2q})$, with $4p + 2q = \sum n_i$, and the subspace $\Pi_\lambda^{(n_i)}$ has dimension over $F$ equal to the number of d*-tableaux of content $(n_i)$ and shape $\lambda$.*

The multilinear component of total degree $n$ of the graded algebra $B(I_w)$ is usually denoted as $\Gamma_n(I_w)$. Since $c_3 \in I_w$ we know that $\Gamma_n(I_w) = 0$ if $n$ is odd.

**Lemma 5.2** *An isomorphism $\varphi$ is given between the vector space $\Gamma_{2m}(I_w)$ and the multilinear component $\Pi^{(1,\ldots,1)}$ of total degree $2m$ of the graded algebra $\Pi$.*

**Proposition 5.3** *The set of proper polynomials $f_S$, as $S$ ranges in the set $N = N^{(1,\ldots,1)}$ of all multilinear normal c-arrays, is linearly independent modulo $I_w$.*

*Proof:* By Theorem 3.3 ($I \subset I_w$), we have that the polynomials $f_S$, for all $S \in N$, form a generating set of the space $\Gamma_{2m}(I_w)$. Hence, their images $\varphi(f_S)$ under the isomorphism $\varphi$ linearly span $\Pi^{(1,\ldots,1)}$. By Propositions 4.6 and 4.7, the number of such polynomials equals the number of multilinear d-tableaux of any shape $(2^{2p}, 1^{2q})$, with $4p + 2q = 2m$. Owing to Proposition 5.1, such number is just the dimension of $\Pi^{(1,\ldots,1)} \approx \Gamma_{2m}(I_w)$. It follows that the elements $f_S$, for all $S \in N$, are linearly independent modulo $I_w$. ∎

For proving Lemma 5.2, it is sufficient to observe what follows. Recall that $R = M_{1,1}(E)$ is a $\mathbb{Z}_2$-graded algebra with homogeneous components:

$$R_0 = \{\begin{pmatrix} a & 0 \\ 0 & d \end{pmatrix} \mid a, d \in E_0\} \quad \text{and} \quad R_1 = \{\begin{pmatrix} 0 & b \\ c & 0 \end{pmatrix} \mid b, c \in E_1\}$$

Note that $W = (W \cap R_0) \oplus (W \cap R_1)$. Since $W \cap R_0$ is in the center of $R$, a proper multilinear polynomial $f(x_1, \ldots, x_n) \in F\langle X \rangle$ is a weak identity for $R$ if and only if it vanishes for all substitutions of indeterminates with elements in $W \cap R_1 = R_1$. We may consider $f(z_1, \ldots, z_n)$ as an element of the free superalgebra $F\langle Y \cup Z \rangle$, where $Y$ is the set of even variables and $Z$ the set of odd ones. We have that $f(x_1, \ldots, x_n)$ is a weak identity for $R$ if and only if $f(z_1, \ldots, z_n)$ is a $\mathbb{Z}_2$-graded polynomial identity for the superalgebra $R$. For a fixed superalgebra $A$, the ideal of $F\langle Y \cup Z \rangle$ of the $\mathbb{Z}_2$-graded polynomial identities of $A$ is usually denoted as $T_2(A)$. Let $V_{p,q}$ denote the subspace of the algebra $F\langle Y \cup Z \rangle$ given by all multilinear polynomials in the variables $y_1, \ldots, y_p$ and $z_1, \ldots, z_q$, and define $\Gamma_{p,q} \subset V_{p,q}$ the subspace of the proper multilinear polynomials. By putting $\Gamma_{p,q}(T_2(A)) = \Gamma_{p,q}/(\Gamma_{p,q} \cap T_2(A))$, for any even integer $2m$ we have:

**Remark 5.4** An isomorphism between the vector space $\Gamma_{2m}(I_w)$ and the space $\Gamma_{0,2m}(T_2(R))$ is simply defined by $x_i \mapsto z_i$.

Consider now $M_2 = M_2(F)$ endowed with its natural $\mathbb{Z}_2$-graduation:

$$(M_2)_0 = \{ \begin{pmatrix} a & 0 \\ 0 & d \end{pmatrix} \mid a, d \in F \} \quad \text{and} \quad (M_2)_1 = \{ \begin{pmatrix} 0 & b \\ c & 0 \end{pmatrix} \mid b, c \in F \}$$

Then, we have that $R = M_2 \hat{\otimes} E = ((M_2)_0 \otimes E_0) \oplus ((M_2)_1 \otimes E_1)$ is the Grassmann envelope of the superalgebra $M_2$. In [13] Kemer describes, for the characteristic zero, an existing relationship between the $\mathbb{Z}_2$-graded polynomial identities of any superalgebra and those of its Grassmann envelope. For the multilinear case, this relationship remains valid in positive characteristic. More precisely, an automorphism "$*$" is defined on the vector space $V_{p,q}$ in the following way. Let $m = u_0 z_{\sigma(1)} u_1 z_{\sigma(2)} u_2 \cdots u_{q-1} z_{\sigma(q)} u_q$ be any monomial of $V_{p,q}$, where $u_0, u_1, \ldots, u_q$ are monomials in the indeterminates $y_i$ and $\sigma \in \mathbb{S}_q$ a suitable permutation. Then, we define:

$$m^* = (-1)^\sigma m$$

**Proposition 5.5** ([13], Lemma 4) Let $A$ be any superalgebra and $f \in V_{p,q}$ a polynomial. It holds:

$$f \in T_2(A) \Leftrightarrow f^* \in T_2(A \hat{\otimes} E)$$

By putting $\Delta_{2m} = \Gamma_{0,2m}^*$ and $\Delta_{2m}(T_2(M_2)) = \Delta_{2m}/\Delta_{2m} \cap T_2(M_2)$, from the previous proposition it follows immediately:

**Remark 5.6** $\Gamma_{0,2m}(T_2(R)) \approx \Delta_{2m}(T_2(M_2))$.

Consider now the polynomial ring $P' = F[H_i, K_i, U_i, V_i]$ in the commuting variables $H_i, K_i, U_i, V_i$, where $i \geq 1$. Note that the relatively free superalgebra $F\langle Y \cup Z \rangle / T_2(M_2)$ is canonically isomorphic to the subalgebra of $M_2(P')$ generated by the matrices:

$$Y_i = \begin{pmatrix} H_i & 0 \\ 0 & K_i \end{pmatrix} \quad \text{and} \quad Z_i = \begin{pmatrix} 0 & U_i \\ V_i & 0 \end{pmatrix}$$

Under the sequence of isomorphisms we defined, if $S$ is a multilinear array of type:

$$S = \begin{bmatrix} a_1 & a_2 & \ldots & a_m \\ b_1 & b_2 & \ldots & b_m \end{bmatrix}$$

each element $f_S + I_w \in \Gamma_{2m}(I_w)$ maps to the scalar matrix $(-1)^S p_S \cdot I_2$, where $(-1)^S$ is the sign of the permutation:

$$\begin{pmatrix} 1 & 2 & \ldots & 2m-1 & 2m \\ a_1 & b_1 & \ldots & a_m & b_m \end{pmatrix}$$

and $p_S$ is the polynomial of $P = F[U_i, V_i \mid 1 \leq i \leq 2m] \subset P'$ defined as:

$$p_S = (U_{a_1} V_{b_1} + U_{b_1} V_{a_1}) \cdots (U_{a_m} V_{b_m} + U_{b_m} V_{a_m})$$

Note that, under the following automorphism of the ring $P$:

$$U_i \mapsto \frac{U_i + \epsilon V_i}{\sqrt{2}}, V_i \mapsto \frac{U_i - \epsilon V_i}{\sqrt{2}}, \quad \text{where } \epsilon^2 = -1$$

the polynomial $p_S$ maps finally to:

$$q_S = q_{a_1 b_1} \cdots q_{a_m b_m}$$

We conclude that the claimed isomorphism $\varphi$ of Lemma 5.2 is defined by putting $\varphi(f_S + I_w) = (-1)^S q_S$.

## 6 Reduction to the multilinear case

**Proposition 6.1** *A one-to-one correspondence is given between the normal c-arrays of the set $N = N^{(n_1, \ldots, n_k)}$ and the ones of the set $N' = N^{(n_{\sigma(1)}, \ldots, n_{\sigma(k)})}$, for any $\sigma \in \mathbb{S}_k$.*

*Proof:* By Propositions 4.6 and 4.7, we have bijection between $N$ and the set of all the d-tableaux of content $(n_1, \ldots, n_k)$ and shape $\lambda = (2^{2p}, 1^{2q})$, for any $p, q$. The same happens to $N'$. Then, the claim follows from an equivalent result for the semistandard tableaux ([11], Proposition 4.3.2). ■

**Proposition 6.2** *Let $l > 0$ denote the number of times the integer 1 occurs in the content $(n_i) = (n_1, \ldots, n_k)$ and let $q$ be the number of occurrences of 2. If $q$ is even then we have a bijection between $N^{(n_i)}$ and $N^{(1, \ldots, 1)}$ (the number of 1 in the content $(1, \ldots, 1)$ is $l$). Otherwise, the bijection is given between $N^{(n_i)}$ and $N^{(1, \ldots, 1, 2)}$.*

*Proof:* Let $q$ be an even integer. By the previous proposition, we may assume that:

$$(n_1, \ldots, n_k) = (\underbrace{2, \ldots, 2}_{r}, \underbrace{1, \ldots, 1}_{2s}, \underbrace{2, \ldots, 2}_{r})$$

where $2r = q$, $2s = l$ and therefore $k = 2r + 2s$. Let $S$ be any c-array of $N^{(n_i)}$. It is sufficient to note that, owing to conditions ($s_1$)–($s_4$), the normal c-array $S$ is necessarily of type:

$$S = \begin{bmatrix} a_1 & \ldots & a_s & k-r+1 & k-r+1 & \ldots & k & k \\ b_1 & \ldots & b_s & r & r & \ldots & 1 & 1 \end{bmatrix}$$

where $\{a_1, b_1, \ldots, a_s, b_s\} = \{r+1, \ldots, r+2s\}$, and the array:

$$S' = \begin{bmatrix} a_1 & \ldots & a_s \\ b_1 & \ldots & b_s \end{bmatrix}$$

is an element of $N^{(0,\ldots,0,1,\ldots,1,0,\ldots,0)}$. In fact, by condition $(s_1),(s_2)$ the integer $k$ has to occur in the last two entries of the first row of $S$. The integer 1 is necessarily in the second row, and by conditions $(s_2),(s_4)$, it has to occur in the last two entries of the second row, and so on. The claimed bijection is the one trivially deduced from $S \mapsto S'$. In the same way, we argue for $q$ odd. ∎

**Proposition 6.3** *Let $l > 0$ be the number of occurrences of 1 in the content $(1,\ldots,1,2)$. A bijection is given between the sets $N = N^{(1,\ldots,1,2)}$ and $N^{(1,\ldots,1)}$.*

*Proof:* Any c-array $S \in N$ is of type:

$$S = \begin{bmatrix} a_1 & \ldots & a_s & l+1 & l+1 \\ b_1 & \ldots & b_s & x & y \end{bmatrix}$$

where $x \leq y$ by condition $(s_2)$. From $(s_1)$ it follows that the integer 1 is in the second row and hence $x = 1$ (otherwise the condition $(s_4)$ is violated). By putting:

$$S' = \begin{bmatrix} a_1 & \ldots & a_s & l+1 \\ b_1 & \ldots & b_s & y \end{bmatrix}$$

we have that $S' \in N^{(0,1,\ldots,1)}$ and the bijection is the one trivially deduced from $S \mapsto S'$. ∎

Note finally that the number of elements of the set $N^{(2,\ldots,2)}$ is just 1 or 0 whenever the number of 2 is respectively even or odd. For concluding the argument of the linear independence we need the following results.

**Lemma 6.4** *The set of proper polynomials $f_S$, as $S$ ranges in the set $N = N^{(1,\ldots,1,2)}$, is linearly independent modulo $I_w$.*

*Proof:* By the proof of Proposition 6.3, we have that any c-array $S \in N$ is of type:

$$S = \begin{bmatrix} a_1 & \ldots & a_s & l+1 & l+1 \\ b_1 & \ldots & b_s & 1 & y \end{bmatrix}$$

Then, the linearization of the polynomial $f_S$ is given modulo $I_w$ by the sum of the polynomials $f'_S, f''_S$ corresponding to the following normal c-arrays:

$$S' = \begin{bmatrix} a_1 & \ldots & a_s & l+1 & l+2 \\ b_1 & \ldots & b_s & 1 & y \end{bmatrix} \qquad S'' = \begin{bmatrix} a_1 & \ldots & a_s & l+1 & l+2 \\ b_1 & \ldots & b_s & y & 1 \end{bmatrix}$$

By the relation $\sum_i \alpha_i f_{S_i} = 0 \mod I_w$, we get hence:

$$\sum_i \alpha_i f_{S'_i} + \sum_i \alpha_i f_{S''_i} = 0$$

Since the polynomials $f_{S'_i}, f_{S''_i}$ are all multilinear and distinct, from Proposition 5.3 it follows that $\alpha_i = 0$ for any $i$. ∎

**Lemma 6.5** Let $\bar{S}$ be a normal c-array of content $(n_i) = (n_1, \ldots, n_k)$, where $n_i = 0$ or $n_i = 2$. Denote by $f_{\bar{S}}$ the proper polynomial associated to $\bar{S}$ and let $f$ be any polynomial of $F\langle X \rangle$ whose variables are disjoint from those of $f_{\bar{S}}$. It holds: $f \cdot f_{\bar{S}} \in I_w$ if and only if $f \in I_w$.

*Proof:* We may assume that:

$$\bar{S} = \begin{bmatrix} r+1 & r+1 & \ldots & 2r & 2r \\ r & r & \ldots & 1 & 1 \end{bmatrix}$$

and $f = f(x_{2r+1}, \ldots, x_n)$ with $r, n > 0$ integers. For any substitution of the variables $x_{2r+1}, \ldots, x_n$ with matrices $w_{2r+1}, \ldots, w_n \in W$, we denote:

$$f(w_{2r+1}, \ldots, w_n) = \begin{pmatrix} a & b \\ c & d \end{pmatrix} \qquad (a, b, c, d \in E)$$

Since the Grassmann algebra $E$ is generated by a vector space $V$ of countable dimension, we may found $4r$ vectors $u_1, \ldots, u_{2r}, v_1, \ldots, v_{2r}$ of the basis of $V$ which do not occur in $a, b, c, d$. Then, we put:

$$w_i = \begin{pmatrix} 0 & u_i \\ v_i & 0 \end{pmatrix} \qquad (i = 1, 2, \ldots, 2r)$$

An easy calculation shows that $f_{\bar{S}}(w_1, \ldots, w_{2r})$ is the scalar matrix:

$$2^r u_1 \cdots u_{2r} v_1 \cdots v_{2r} \cdot I_2$$

Hence, from $f(w_{2r+1}, \ldots, w_n) \cdot f_{\bar{S}}(w_1, \ldots, w_{2r}) = 0$ it follows necessarily that $f(w_{2r+1}, \ldots, w_n) = 0$. ∎

## 7 Main results

**Theorem 7.1** Denote by $l$ and $q$ respectively the number of times the integers 1 and 2 occurs in the content $(n_i) = (n_1, \ldots, n_k)$, and put $n = \sum n_i$. The multihomogeneous component $B^{(n_i)}(I_w) = 0$ if and only if one of the following conditions is verified:

i) $n$ is odd;

ii) $n$ is even, $n_i > 2$ for some $i$;

iii) $n$ is even, $n_i \leq 2$ for all $i$, $l = 0$, $q$ is odd.

Otherwise, the set of the proper polynomials $f_S$, as $S$ ranges in the set $N^{(n_i)}$, is a linear basis of the space $B^{(n_i)}(I_w)$ whose dimension is $\binom{2s-1}{s}$, where $l = 2s$.

*Proof:* In section 3 we have shown that $B^{(n_i)}(I_w) = 0$ for the cases (i) and (ii). In the remaining case, note that the polynomials $f_S$, for all $S \in N^{(n_i)}$, span the space $B^{(n_i)}(I_w)$ by Theorem 3.3 ($I \subset I_w$). The first part of the claim is then proved since $N^{(n_i)} = \emptyset$ if and only if condition (iii) is true.

If $B^{(n_i)}(I_w) \neq 0$, we have to prove now that the polynomials $f_S$ are linearly independent modulo $I_w$. Note that the vector spaces $B^{(n_1,\ldots,n_k)}(I_w)$ and $B^{(n_{\sigma(1)},\ldots,n_{\sigma(k)})}(I_w)$ are canonically isomorphic for any $\sigma \in \mathbb{S}_k$. On the other hand, by Proposition 6.1 we have a one-to-one correspondence between the sets $N^{(n_1,\ldots,n_k)}$ and $N^{(n_{\sigma(1)},\ldots,n_{\sigma(k)})}$. Then, it is sufficient to prove the linear independence for the case:

$$(n_1,\ldots,n_k) = (\underbrace{2,\ldots,2}_{r},\underbrace{1,\ldots,1}_{2s},\underbrace{2,\ldots,2}_{t})$$

where $t = r$ or $t = r + 1$ when the number $q$ of occurrences of 2 in the content $(n_1,\ldots,n_k)$ is respectively even or odd. As in the proof of Proposition 6.2, it holds that any c-array $S \in N^{(n_i)}$ is of type $S = S' * \bar{S}$, where:

$$\bar{S} = \begin{bmatrix} k-r+1 & k-r+1 & \ldots & k & k \\ r & r & \ldots & 1 & 1 \end{bmatrix} \quad \text{and} \quad S' = \begin{bmatrix} a_1 & \ldots & a_s \\ b_1 & \ldots & b_s \end{bmatrix}$$

The c-array $S'$ belongs to $N^{(0,\ldots,0,1,\ldots,1,0,\ldots,0)}$ or $N^{(0,\ldots,0,1,\ldots,1,2,0,\ldots,0)}$ whenever $q$ is even or odd. Note that in both cases the sets of the entries of $\bar{S}$ and $S'$ are disjoint and we have $f_S = f'_S \cdot f_{\bar{S}}$. Suppose now that $\sum_i \alpha_i f_{S_i} = 0$ mod $I_w$. Then $(\sum_i \alpha_i f_{S'_i}) \cdot f_{\bar{S}} = 0$ and hence $\sum_i \alpha_i f_{S'_i} = 0$ owing to Lemma 6.5. If the integer $q$ is even, the polynomials $f_{S'_i}$ are all multilinear and linearly independent modulo $I_w$ by Proposition 5.3. Otherwise, the claim follows by Lemma 6.4.

For the case $(n_i) = (2,\ldots,2)$, where $q$ even, we have immediately that the dimension $\dim_F B^{(n_i)}(I_w) = |N^{(n_i)}| = 1$. In the remaining cases in which $B^{(n_i)}(I_w) \neq 0$, that is for $l > 0$, by Propositions 6.2 and 6.3 we have that the normal c-arrays of $N^{(n_i)}$ are in ono-to-one correspondence with the multilinear normal c-arrays of $N^{(1,\ldots,1)}$, (the number of 1 in the content $(1,\ldots,1)$ is $l$). Then, the dimension $\binom{2s-1}{s}$ can be easily computed by enumerating all the multilinear d-tableaux of shape $(2^{2p}, 1^{2q})$, with $4p + 2q = 2s$, that correspond under the algorithm `Carray2Dtableau` to the elements of $N^{(1,\ldots,1)}$. For an alternative approach based on a straightforward enumeration of the multilinear normal c-arrays we refer to [5]. ∎

**Theorem 7.2** *Let $F$ be an infinite field of characteristic different from 2. The ideal $I_w$ of the weak polynomial identities of the superalgebra $M_{1,1}(E)$ is $\Omega$-generated by the polynomials:*

$$c_3 = [x_1, x_2, x_3] \quad \text{and} \quad p = [x_2, x_1][x_3, x_1][x_4, x_1]$$

*Moreover, the Hilbert series of the algebra $B_k(I_w)$ is the following:*

$$H(B_k(I_w), t_1, \ldots, t_k) = \frac{1}{2} \sum_{i=0}^{k} [e_i^2(t_1, \ldots, t_k) + (-1)^i e_i(t_1^2, \ldots, t_k^2)]$$

*where $e_i(t_1, \ldots, t_k)$ is the $i$-th elementary symmetric function.*

*Proof:* Recall that $I \subset I_w$ is the ideal $\Omega$-generated by the polynomials $c_3, p$. By means of Theorems 3.3 and 7.1, for any multidegree $(n_i) = (n_1, \ldots, n_k)$ we have shown that $B^{(n_i)}(I) = B^{(n_i)}(I_w)$ and therefore $I = I_w$.

Owing to Propositions 4.6 and 4.7, the basis of the vector space $B_k^{(n_i)}(I_w)$ described in the Theorem 7.1 implies the following decomposition:

$$B_k^{(n_i)}(I_w) = \bigoplus_\lambda B_{k,\lambda}^{(n_i)}(I_w)$$

where $\lambda$ ranges in the set (eventually empty) of the partitions of type $\lambda = (2^{2p}, 1^{2q})$, with $4p + 2q = \sum n_i$, and the subspace $B_{k,\lambda}^{(n_i)}$ has dimension over $F$ equal to the number of d-tableaux of content $(n_i)$ and shape $\lambda$. Let $s_\lambda(t_1, \ldots, t_k)$ be the Schur function corresponding to the shape $\lambda$. By the previous decomposition it holds:

$$H(B_k(I_w), t_1, \ldots, t_k) = \sum_\lambda s_\lambda(t_1, \ldots, t_k)$$

where $\lambda = (2^{2p}, 1^{2q})$, for all the integers $p, q \geq 0$. Now, the following formula is due to Carini and Drensky ([2], page 478):

$$\sum_\lambda s_\lambda(t_1, \ldots, t_k) = \frac{1}{2} \sum_{i=0}^{k} [e_i^2(t_1, \ldots, t_k) + (-1)^i e_i(t_1^2, \ldots, t_k^2)]$$

∎

For the sake of completeness, we finally add a result which has been proved in [5] (Theorem 5.4).

**Proposition 7.3** *The generating function $\gamma(I_w, z)$ of the proper codimension sequence $\dim_F \Gamma_k(I_w) = \dim_F B_k^{(1,\ldots,1)}(I_w)$, for all $k \geq 0$, is the following:*

$$\gamma(I_w, z) = \frac{1}{2}(1 + \frac{i}{\sqrt{(2z+1)(2z-1)}})$$

We are pleased to thank Vesselin Drensky for the many stimulating discussions during his visit to the University of Bari.